\documentclass[12pt]{amsart}

\usepackage{amsfonts, amssymb, amscd}

\usepackage{verbatim}
\usepackage{amssymb}
\usepackage{mathrsfs}
\usepackage{graphicx}
\usepackage{cite}
\usepackage[all]{xy}
\usepackage{tikz}
\usepackage{MnSymbol}

\usepackage{hyperref}

\newcommand{\Zz}{\mathbb{Z}}

\newcommand{\Pp}{\mathbb{P}}
\newcommand{\Rr}{\mathbb{R}}
\newcommand{\Qq}{\mathbb{Q}}

\newcommand{\Exc}{\operatorname{Exc}}
\newcommand{\Mov}{\operatorname{Mov}}
\newcommand{\Bs}{\operatorname{Bs}}

\newcommand{\Oo}{\mathcal{O}}

\newcommand\blfootnote[1]{%
  \begingroup
  \renewcommand\thefootnote{}\footnote{#1}%
  \addtocounter{footnote}{-1}%
  \endgroup
}

\numberwithin{equation}{section}

\begin{document}

\title[Movable cones]{Counterexamples of Lefschetz hyperplane type results for movable cones}

\blfootnote{ {\it Acknowledgements}: The author benefits from many
conversations with Professors Lev Borisov, Chenyang Xu and John
Christian Ottem. Thanks also go to Professor Yoshinori Gongyo who
pointed out the reference \cite{CO15}.}

\begin{abstract}
The main theorem of the paper provides a way to produce examples such that the movable cone of an ample divisor does not coincide with the movable cone of its ambient variety.
\end{abstract}

\author{Zhan Li}
\address{Beijing International Center for Mathematical Research\\
Beijing 100871, China} \email{lizhan@math.pku.edu.cn}

\maketitle

The movable cone of a variety $X$ is defined to be the convex cone in $H^2(X, \Rr)$ generated by the classes of movable divisors. Structures of a movable cone carry information of birational classes. For example, the Morrison-Kawamata cone conjecture roughly states that there is a fundamental domain inside the intersection of movable and effective cones for the birational action of a Calabi-Yau variety (see \cite{Kaw97} for the precise statement and partial results). Moreover, given many successful induction arguments (adjunctions, pl-flips, special terminations etc.) in the minimal model program involving the investigations of properties of a divisor and its ambient space, we wish to know whether their movable cones are compatible in any sense. One version can be obtained by analogizing the following Lefschetz hyperplane theorem:

\medskip
Suppose $P$ is a smooth projective variety with $\dim P \geq 4$, and $X \subseteq P$ is an ample divisor. Then the natural map $H^p(P, \Qq) \to H^p(X, \Qq)$ induced by the inclusion is an isomorphism  for $p \leq \dim P - 2$ and injective for $p = \dim P-1$.
\medskip

We call the analogy of this theorem for other structures the Lefschetz hyperplane type results. Lefschetz hyperplane type results hold for fundamental groups, Class groups with general ample divisors (\cite{RS06}) and even stacks (\cite{HL10}), but fail for ample cones (\cite{HLW02}, see Remark 2\ref{remark} for some positive results). Theorem 1\ref{main theorem} given below manifests that one cannot expect Lefschetz hyperplane type result holds for movable cones even for general ample divisors.

\medskip

Before stating the theorem, let us fix the notations and terminologies:

The movable divisor $D$ is a Cartier divisor whose complete linear system $|D| \neq 0$ and the base locus $\Bs(D)$ has codimension bigger than $1$. In the same fashion, one can define a movable divisor in the relative setting (\cite{Kaw88} \S 2). Let $\Mov(X) \subseteq H^2(X, \Rr)$ be the movable cone of $X$, that is, the rational convex cone generated by the movable divisor classes $[D]$. Let $\overline \Mov(X)$ and $\Mov(X)^\circ$ be the closure and the interior of the cone $\Mov(X)$ respectively. Suppose $g: X \to Y$ is a morphism between normal varieties with connected fibers. Then either $\dim X > \dim Y$ or $f$ is a birational morphism. To coin a term from minimal model program, we call $f$ a \emph{fibre morphism} if $\dim X > \dim Y$; a \emph{divisorial morphism} if $f$ is birational and the exceptional locus $\Exc(f)$ of $f$ is of dimensional $n-1$; and a \emph{small morphism} if $f$ is birational with $\dim \Exc(f) \leq \dim X - 2$.

\medskip

{\noindent\bf Theorem 1.\label{main theorem}}
{\it Let $P$ be a normal, projective variety of dimension $n$ and $X \subseteq P$ be a normal Weil divisor. Suppose there exists a surjective, projective morphism $p: P \to Q$ with connected fibers. Let $f: X \to Z$ be the Stein factorization of the restriction morphism $p|_X$.

(1) If $p$ is a fibre or divisorial morphism, and $f$ is a small morphism (see the explanations above), then there exists a Cartier divisor $S$ on $P$ such that $[S] \notin \overline{\Mov}(P)$, but $[S|_X] \in \Mov(X)^\circ$.

(2) If $p$ is a small morphism, and $f$ is a divisorial morphism, then there exists a Cartier divisor $S$ on $P$ such that $[S] \in \Mov(P)^\circ$, but $[S|_X] \notin \overline \Mov(P)$.
}

\begin{proof}
The morphisms stated in the theorem are labeled on the following diagram
\[
\xymatrix{
P \ar[rrd]^p \\
X \ar@{^{(}->}[u] \ar[rrd]_{p|_X}  \ar@/^/[r]^f& Z \ar@/^/[rd] & Q\\
&& p(X) \qquad.
}
\]

For (1), let $H$ be an ample divisor on $Q$ such that $p: P \to Q$ is defined by the base point free linear system $|p^*H|$. Then the morphism $f: X \to Z$ is the Stein factorization of the morphism defined by the linear system $|(p^*H)|_X|$. Let $S'$ be an ample divisor on $P$.

\medskip

We first claim that $[(p^*H)|_X]$ lives in the interior of $\Mov(X)$. In fact, let $A_Z$ be an ample divisor on $Z$ such that $f^*A_Z = (p^*H)|_X$. For any ample divisor $B_X$ on $X$, let $B_Z=f_*B_X$ be its strictly transform on $Z$. There exists $m \gg 0$, such that $\Oo_Z(-B_Z)(mA_Z)$ is generated by its global sections, hence base point free. The base locus of $\Oo_Z(-B_X)(m(p^*H)|_X)$ can only be contained in the exceptional locus of $f$. Because $f$ is a small contraction, $\Bs(\Oo_Z(-B_X)(m(p^*H)|_X))$ is at most of codimensional $2$. This shows that $m(p^*H)|_X - B_X$ is a movable divisor, that is $[(p^*H)|_X - \frac{1}{m} B_X] \in \Mov(X)$. As a result $[(p^*H)|_X] \in \Mov(X)^\circ$.

\medskip

Because $[(p^*H)|_X] \in \Mov(X)^\circ$, there exists $N$ such that $n > N$, $[n(p^*H)|_X - S'|_X] \in \Mov(X)^\circ$. Let $S = n(p^*H) - S'$, we will show that for $n$ sufficiently large, $[S] \notin \overline\Mov(P)$ and thus complete the proof.

\medskip

Suppose otherwise, $[S] \in \overline\Mov(P)$. For any $m \geq n$, we have $m(p^*H) - S = (m-n) p^*H  + S'$ to be an ample divisor on $P$, and hence $[m(p^*H) - S] \in \Mov(P)^\circ$. Because we assume $[S] \in \overline\Mov(P)$, $[m(p^*H) - S] + [S] \in \Mov(P)^\circ$. For any ample divisor $\Theta$ on $P$, there exists $0<\delta \ll1$ such that
\[
[m(p^*H) -\delta \Theta]=[m(p^*H) - S] + [S] + [-\delta \Theta] \in \Mov(P).
\] Hence, there exists sufficiently divisible $l>0$ such that $l \left(m(p^*H) -\delta \Theta\right)$ is a movable divisor. However, any curve $C$ contracted by $p$ has intersection
\[l \left(m(p^*H) -\delta \Theta\right) \cdot C = -l \delta \Theta \cdot  C < 0.
\] Hence $C \subseteq \Bs(l \left(m(p^*H) -\delta \Theta\right))$. Because $\Exc(p)$ is covered by curves contracted by $p$, we have $\Exc(p) \subseteq \Bs(l \left(m(p^*H) -\delta \Theta\right))$. This is a contradiction since $\dim (\Exc(p)) \geq \dim P -1$ but the dimension of $ \Bs(l \left(m(p^*H) -\delta \Theta\right))$ is at most $ \dim P - 2$.

\medskip

The claim (2) can be proved similarly.  We just sketch the argument. Let $H, S', S$ and $A_Z$ be chosen as before. Because $p$ is assumed to be small morphism, $[n p^*H - S'] \in \Mov(P)^\circ$ for $n \gg 0$. However, $S|_X = (n p^*H - S')|_X$ cannot correspond to a class live in $\overline\Mov(X)$ for sufficiently large $n$. In fact, because $S'|_X$ is ample, for $m \geq n$,  $m (p^*H)|_X - S|_X = (m-n)(p^*H)|_X + S'|_X$ is also an ample divisor. As a result, $m (p^*H)|_X - S|_X \in \Mov(X)^\circ$. If $[S|_X] \in \overline\Mov(X)$, there exists $0< \xi \ll 1$ and an ample divisor $\Xi$ on $X$ such that
\[
[m (p^*H)|_X - \xi \Xi] =[m (p^*H)|_X - S|_X] + [S|_X] + [-\xi\Xi] \in \Mov(X).
\] However, this will give a contradiction since there exists $l$ such that $l\left(m (p^*H)|_X - \xi \Xi\right)$ is movable but $\Exc(f) \subseteq \Bs(l\left(m (p^*H)|_X - \xi \Xi\right))$ is of dimension $\dim X -1$.
\end{proof}

{\noindent \bf Remark 2.}\label{remark}
{\it This theorem also holds in the relative setting without any change of argument. That is, one can consider the morphisms over a scheme, and use relative movable divisors in the place of movable divisors.
}

\medskip

Given this result, it is easy to construct counterexamples of Lefschetz hyperplane type result for movable cones. In fact, such result does not hold even for generic ample divisors.

\medskip

{\noindent \bf Proposition 3.}\label{proposition}
{\it There are examples for each $n, n\geq 4$ which consist a smooth projective variety $P$ of dimensional $n$, and any generic ample divisor $X$ of $P$, such that $\overline\Mov(X)$ and $\overline\Mov(P)$ do not coincide under the isomorphism $H^2(P, \Rr) \cong H^2(X, \Rr)$.
}

\begin{proof}
Suppose $p: P \to P'$ is a blowup of a smooth codimensional $2$ subvariety $W$ in a smooth projective $n$ dimensional  variety $P'$. Then $P$ is a smooth projective variety and the fibre of $p^{-1}(W) \to W$ is $\Pp^1$. By Bertini theorem, any generic ample divisor $X$ is smooth and intersects the general fibre of $p^{-1}(W) \to W$ at finite points. Then the Stein factorization $f: X \to Z$ of $p|_X: X \to p(X)$ is a small morphism. In fact, by the choice of $X$, the exceptional locus of $f$ is contained in $\{x \in X \mid \dim p^{-1}(p(x)) \geq 1\}$, hence at most of dimension $\dim (p^{-1}(W) \cap X) -1= n-3$. Then by Theorem 1\ref{main theorem}, there exists a divisor $S$, such that $ [S] \notin\overline \Mov(P)$ but its restriction $[S|_X]$ lives in the interior of the movable cone of $X$.
\end{proof}

On the other hand, a general $X \subseteq P$ intersects base locus of a divisor $S$ transversally. In this case,  if $S$ is movable then $S|_X$ is also movable because $\Bs(S|_X) \subseteq \Bs(S) \cap X$. This phenomenon is well reflected by Theorem 1\ref{main theorem}(2), that is, $X$ has to contain the exceptional locus of $p$, and hence cannot be general.

\medskip

Koll\'ar showed (see \cite{Kol91}) that for any smooth Fano variety $P$ of dimension greater than $3$, and a divisor $X \subset P$, the natural inclusion of numerically effective cones $i_*: {\rm NE}(P) \to {\rm NE}(X)$ is an isomorphism. However, even in this case, $\Mov(P)$ and $\Mov(X)$ could still be different: there are extremal contractions of Fano manifold whose general fibre of exceptional divisor is of dimension $1$. Then the previous construction will give non-movable divisor on $P$ whose restriction to any generic ample divisor (in particular generic $X \in |-K_P|$) is movable.

\medskip

John Ottem pointed out that a simple example of the same kind can be
obtained by considering hypersurfaces in the product of projective
spaces. For example, let $X$ be a general bidegree $(2,1)$
hypersurface in $\Pp^1 \times \Pp^3$. Written in homogenous
coordinates, it is defined by $x^2_0 f_0 + x_0 x_1 f_1 + x_1^2 f_2
=0$. The Picard group of $X$ is isomorphic to $\Zz^2$ by Lefschetz
hyperplane theorem. The second projection $\rm pr_2: X \to \Pp^3$ is
a double cover outside of $\{f_0 = f_1 = f_2 =0\} \subseteq \Pp^3$.
Let $\sigma: X \dashedrightarrow X$ be the map defined by switching
two sheets. To be precise, it sends
\[
[x_0: x_1 \mid y_0: y_1: y_2] \to [\frac{f_2}{x_0}: \frac{f_0}{x_1}
\mid y_0: y_1: y_2].
\] This is a well defined map outside of a curve $\Pp^1 \times \{f_0 = f_1 = f_2
=0\}$ (we interchange $-(\frac{f_1}{x_0} + \frac{x_1 f_2}{x_0^2})$
and $\frac{f_0}{x_1}$, etc. when $x_1$ or $x_0$ is $0$). Moreover,
it is a small birational morphism. Let $H_1, H_2$ be the restriction
of $(1,0)$ and $(0,1)$ hypersurfaces to $X$. Then the strict
transform ${\sigma}^{-1}_*H_1$ is linearly equivalent to $H_2 -H_1$.
As restrictions of base point free divisors, $H_1, H_2$ are also
base point free. The strict transform ${\sigma}^{-1}_*H_1 = H_2
-H_1$ is movable because $\sigma$ is small. On the other hand,
$H^0(\Pp^1 \times \Pp^3, \Oo(-m, 2m))=0$ for any $m>0$. In
particular, $[\Oo(-1, 2)] \notin \Mov(\Pp^1 \times \Pp^3)$. In fact,
one can determine the movable cone of $X$ explicitly: because $H_1$
is not a big divisor, ${\sigma}^{-1}_*H_1$ is not a big divisor
either. Since $\Mov^\circ(X)$ consists of big divisors, $H_1$ and
${\sigma}^{-1}_*H_1$ form the two rays generated $\Mov(X)$ and thus
$\overline \Mov(X)$. On the other hand, $\Mov(\Pp^1 \times \Pp^3) =
\overline \Mov(\Pp^1 \times \Pp^3)$ are generated by $[\Oo(1,0)]$
and $[\Oo(0,1)]$. As a result, $\Mov(X)$ and $\Mov(\Pp^1 \times
\Pp^3)$ do not coincide under the natural restriction. We recommend
\cite{Ott14} for detailed discussions of related problems for
hypersurfaces in the products of projective spaces.

\medskip

{\noindent \bf Remark 4.~}{\it Yoshinori Gongyo pointed out that a
similar example with non-isomorphic movable cones had already
appeared in \cite{CO15} (see Remark 4.2).}

\bibliographystyle{alpha}
\bibliography{bibfile}

\end{document}